\documentclass[lettersize,journal]{IEEEtran}

\usepackage{inputenc}
\usepackage[style=ieee, backend=biber, isbn=false]{biblatex}

\usepackage{amsmath,amsfonts}
\usepackage{algorithmic}
\usepackage{algorithm}
\usepackage{array}
\usepackage[caption=false,font=normalsize,labelfont=sf,textfont=sf]{subfig}
\usepackage{textcomp}
\usepackage{stfloats}
\usepackage{url}
\usepackage{verbatim}
\usepackage{graphicx}
\usepackage{todonotes}
\usepackage{comment}

\usepackage{amssymb}
\usepackage{latexsym}
\usepackage{sansmath}

\usepackage{url}
\usepackage{xcolor}
\definecolor{newcolor}{rgb}{.8,.349,.1}

\usepackage{comment}

\usepackage{import}

\usepackage[utf8]{luainputenc}

\usepackage[T1]{fontenc}

\usepackage{graphicx}

\usepackage[english]{babel}
\addto\captionsenglish{}
\addto\captionsenglish{}
\usepackage{csquotes}

\usepackage{newtxtext}
\usepackage{amsthm}
\usepackage[slantedGreek]{newtxmath}
\usepackage[OMLmathsfit]{isomath}
\DeclareMathAlphabet{\mathbfsf}{\encodingdefault}{\sfdefault}{bx}{n}
\usepackage{bm}
\usepackage{envmath}
\usepackage{mathtools}
\usepackage{commath}
\usepackage{siunitx}

\usepackage{booktabs}
\usepackage{footmisc}  

\usepackage{url}

\theoremstyle{definition}

\theoremstyle{plain}

\theoremstyle{remark}

\usepackage{lineno}
\modulolinenumbers[5]
\usepackage{umoline}

\usepackage{pgfplots}
\usepackage{pgfplotstable}
\pgfplotsset{compat=newest}
\pgfplotsset{plot coordinates/math parser=false}
\newlength\figureheight
\newlength\figurewidth
\pgfplotsset{every axis plot/.append style={line width=1.5pt},
    legend style={font=\footnotesize, 
        text height=1.0ex,
        draw=black,
        fill=white,
        legend cell align=left}}

\usepackage{hyperref} 
\usepackage[english]{cleveref}

\Crefname{defn}{definition}{definitions}
\Crefname{defn}{Definition}{Definitions}

\Crefname{asm}{assumption}{assumptions}
\Crefname{asm}{Assumption}{Assumptions}

\crefname{lem}{lemma}{lemmas} 
\Crefname{lem}{Lemma}{Lemmas}

\crefname{prop}{proposition}{propositions} 
\Crefname{prop}{Proposition}{Propositions}

\crefname{thm}{theorem}{theorms} 
\Crefname{thm}{Theorem}{Theorms}

\crefname{cor}{corollary}{corollaries}
\Crefname{cor}{Corollary}{Corollaries}
\newcounter{subequation}
\newlength\mtabskip\mtabskip=-1.25cm

\def\mtabLong{long}

\newcommand{\mr}{\mathrm}

\newcommand{\veg}[1]{\bm{#1}}     
\newcommand{\mat}[1]{\mathsfbfit{#1}} 
\renewcommand{\vec}[1]{\mathsfbfit{#1}} 
\newcommand{\op}[1]{\mathcal{#1}} 
\newcommand{\vecop}[1]{\bm{\mathcal{#1}}} 

\newcommand{\ii}{\mathrm{i}}   

\DeclareMathOperator{\ee}{e}

\DeclareMathOperator{\Ai}{Ai}
\DeclareMathOperator{\Bi}{Bi}

\newcommand{\NGl}{\mathrm{}}
\newcommand{\Gl}{\mathrm{G}}

\newcommand{\T}{\mr{T}}
\newcommand{\Ht}{\dagger}

\newcommand\restr[2]{{
        \left.\kern-\nulldelimiterspace 
        #1 
        \vphantom{|} 
        \right|_{#2} 
}}

\newcommand\rst[3]{{
        \left.\kern-\nulldelimiterspace 
        #1 
        \vphantom{|} 
        \right|_{#2}^{#3} 
}}

\usepackage{acro}

\DeclareAcronym{DG}
{
    short = DG ,
    long = discontinuous Galerkin
}

\DeclareAcronym{ACA}
{
    short = ACA ,
    long = adaptive cross approximation
}

\DeclareAcronym{EFIE}
{
    short =  EFIE ,
    long = electric field integral equation
}

\DeclareAcronym{MFIE}
{
    short =  MFIE ,
    long = magnetic field integral equation
}

\DeclareAcronym{CFIE}
{
    short =  CFIE ,
    long = combined field integral equation
}

\DeclareAcronym{MUIE}
{
    short =  MUIE ,
    long = Müller integral equation
}

\DeclareAcronym{PMCHWT}
{
    short =  PMCHWT ,
    long = Poggio-Miller-Chang-Harrington-Wu-Tsai integral equation
}

\DeclareAcronym{SPD}
{
    short =  SPD ,
    long = {symmetric, positive definite}
}

\DeclareAcronym{SPSD}
{
    short =  SPD ,
    long = {symmetric, positive semi-definite}
}

\DeclareAcronym{PEC}
{
    short =  PEC ,
    long = perfectly electrically conducting
}

\DeclareAcronym{RWG}
{
    short = RWG ,
    long = Rao-Wilton-Glisson
} 

\DeclareAcronym{BC}
{
    short = BC ,
    long = Buffa-Christiansen
}

\DeclareAcronym{SVD}
{
    short = SVD ,
    long = singular value decomposition
}

\DeclareAcronym{CG}
{
    short = CG ,
    long = conjugate gradient
} 

\DeclareAcronym{PCG}
{
    short = PCG ,
    long = preconditioned conjugate gradient
} 

\DeclareAcronym{CGS}
{
    short = CGS ,
    long = conjugate gradient squared
}

\DeclareAcronym{CMP}
{
    short = CMP ,
    long = Calderón multiplicative preconditioner
} 

\DeclareAcronym{RFCMP}
{
    short = RF-CMP ,
    long = refinement-free Calderón multiplicative preconditioner
} 

\DeclareAcronym{HPD}
{
    short = HPD ,
    long = {Hermitian, positive definite}
} 

\DeclareAcronym{RHS}
{
    short = RHS ,
    long = {right-hand side}
}

\DeclareAcronym{PW}
{
    short = PW ,
    long = {plane wave}
} 

\DeclareAcronym{HD}
{
    short = HD ,
    long = {Hertzian dipole}
} 

\DeclareAcronym{FF}
{
    short = FF ,
    long = {far-field}
} 

\DeclareAcronym{NF}
{
    short = NF ,
    long = {near-field}
}  

\newcolumntype {n}{c}
\newcolumntype {N}{>{\small}c}
\newcolumntype {L}{>{\small}l}
\newcolumntype {F}{>{\footnotesize}c}
\newcolumntype {v}[1]{>{\raggedright \hspace {0pt}} p {#1}}
\newcolumntype {V}[1]{>{\small \raggedright \hspace {0pt}} p {#1}}
\newcolumntype{d}[1]{>{\DC@{.}{.}{#1}}c<{\DC@end}}

\newcolumntype{R}[1]{%
    >{\begin{turn}{90}\begin{minipage}{#1}\small\raggedright\hspace{0pt}}l%
            <{\end{minipage}\end{turn}}%
}

\NewDocumentCommand{\TA}{o}{
    \IfNoValueTF {#1} {%
        \vecop T_{\kern-2pt\mr{A}}
    }
    {
        \vecop T_{\kern-2pt\mr{A},#1}
    }
}

\NewDocumentCommand{\TPhi}{o}{
    \IfNoValueTF {#1} {%
        \vecop T_{\kern-2pt\Phiup}
    }
    {
        \vecop T_{\kern-2pt\Phiup,#1}
    }
}

\NewDocumentCommand{\matTA}{o}{
    \IfNoValueTF {#1} {%
        \mat T_\mr{A}   
        }
    {
        \mat T_{\mr{A},#1}
    }
}

\NewDocumentCommand{\matTPhi}{o}{
    \IfNoValueTF {#1} {%
        \mat T_\Phiup   
        }
    {
        \mat T_{\Phiup,#1}
    }
}

\NewDocumentCommand{\MSL}{o}{
    \IfNoValueTF {#1} {%
        \veg \Psi_\mr{SL}
        }
    {
        \veg \Psi_{\mr{SL},#1}
    }
}

\NewDocumentCommand{\MDL}{o}{
    \IfNoValueTF {#1} {%
        \veg \Psi_\mr{DL}
        }
    {
        \veg \Psi_{\mr{DL},#1}
    }
}

\NewDocumentCommand{\PA}{o}{
    \IfNoValueTF {#1} {%
        \veg \Psi_\mr{A}
        }
    {
        \veg \Psi_{\mr{A},#1}
    }
}

\NewDocumentCommand{\PPhi}{o}{
    \IfNoValueTF {#1} {%
        \veg \Psi_{\Phiup}
        }
    {
        \veg \Psi_{\Phiup,#1}
    }
}

\begin{document}

\newcommand{\jj}{\mathrm{j}}

\title{\fontsize{24}{28}\selectfont Asymptotic Spectral Insights Behind Fast Direct Solvers for High-Frequency Electromagnetic Integral Equations on Non-Canonical Geometries}


\author{
	\IEEEauthorblockN{Viviana Giunzioni\IEEEauthorrefmark{1}, Clément Henry\IEEEauthorrefmark{2}, Adrien Merlini\IEEEauthorrefmark{2}, and Francesco P. Andriulli\IEEEauthorrefmark{1}}
	
\IEEEauthorblockA{\IEEEauthorrefmark{1} Department of Electronics and Telecommunications, Politecnico di Torino, 10129 Turin, Italy }

\IEEEauthorblockA{\IEEEauthorrefmark{2} Microwave Department, IMT Atlantique, 29238 Brest, France}
	}



\pagenumbering{gobble}

\maketitle

\begin{abstract}   
Integral-equation-based fast direct solvers for electromagnetic scattering can substantially reduce computational costs, especially in the presence of multiple excitations. We recently proposed a new high-frequency fast direct solver strategy that combines preconditioning techniques with acceleration algorithms. However, the validity of this approach applied to non-canonical geometries requires further justification. In this contribution, we collect relevant semiclassical microlocal results and use them to assess the legitimacy and effectiveness of the proposed fast direct solver in the high-frequency regime. 
\end{abstract}


\section{Introduction}
\label{sec:introduction}
\IEEEPARstart{B}{oundary} integral equations discretized via boundary element methods are fundamental tools in the modeling of electromagnetic scattering phenomena. As the frequency of the electromagnetic problem increases, however, the numerical solution becomes increasingly computationally demanding, motivating the research for direct and inverse solution strategies which are ``fast'', that is, methods which reduce the computational time and resources required to achieve the solution with a prescribed accuracy. Fast iterative solvers typically combine preconditioning strategies and acceleration techniques, but require the solution of an entirely new problem whenever the excitation source changes, resulting in a significant inefficiency when the handling of multiple right-hand-sides (RHSs) is needed. In contrast, fast direct solvers aim at constructing a representation of the inverse operator matrix in reduced complexity, making these techniques well suited to this setting.

In this contribution, we aim to provide a justification for a recently introduced filter-based fast direct solver approach for some Calderón combined field integral equations in the high-frequency regime applied to smooth convex 2D boundaries \cite{consoli2022fast, giunzioni2025new}. In particular, the strategy relies on the assumption that the underlying operators, denoted as transverse magnetic and transverse electric Calderón combined field integral operators (TM- and TE-CCFIO), can be written as the sum of a scaled identity and a compact perturbation, denoted as $\op{C}$. The CCFIO, then, can be written as
\begin{equation}
    \text{CCFIO} = \frac{\op I}{2}+\op C  = \frac{\op I}{2} + \op F \left(\text{CCFIO}-\frac{\op I}{2}\right)\,,
    \label{eqn:CCFIOassumption}
\end{equation}
where $\op F$ denotes a spectral filter which retains the full spectral contribution of $\op C$, that is, $\op F (\op C) = \op C$. The independent discretization of the contributions $\op I/2$ and $\op F(\op C)$ allows to correct some elliptic spectral deficiencies which emerge from the discretization of the full operator. The discretization of the compact part can then be approximated in skeleton form with controllable error, resulting, when applying the inverse operator matrix obtained by the Woodbury formula to a RHS, in a solution with controllable error.

Here, to assess the legitimacy of this fast direct solver approach, we show that the spectral content of $\op C$ grows with frequency like $k^{1/3}$, resulting in an increase in computational time of the fast direct solver algorithm at most like $k^{4/3}$, as observed experimentally (see numerical results in \cite{consoli2022fast,giunzioni2025new}). Polarization-dependent compensation phenomena may further reduce the computational burden in the high-frequency regime if a constant solution error is required, so that the above-mentioned increase of $k^{4/3}$ actually constitutes an upper bound.
The strategy adopted to establish these results relies on semiclassical microlocal, stationary-phase, techniques, to infer local spectral properties of the integral operators, as well as local spatial properties of the RHSs and solution functions.

\section{Background and Formalism}
\label{sec:background}

Consider the time-harmonic electromagnetic scattering from a perfect electric conductor (PEC) cylinder indefinitely extended along the longitudinal direction $z$. Let $\Omega$ be the open set modeling the transverse cross-section of the cylinder and $\gamma\coloneq\partial \Omega$ be its 2D convex, smooth, contour. The exterior space $\mathbb{R}^2\backslash\Omega$ is characterized by its impedance $\eta = \sqrt{\mu/\epsilon}$ and the corresponding wavenumber $k = \omega \sqrt{\mu\epsilon}$. We parametrize $\gamma$ along the curvilinear abscissa $s \in [0,L)$ and, for each point along the curve, we identify the outgoing normal and tangential unit vectors to the curve, $n(s)$ and $\tau(s)$, and the curvature $\kappa(s)$.
The single-layer, double-layer, adjoint double-layer, and hypersingular operators are defined as
\begin{align}
    \op S^k f(s) &\coloneqq \int_\gamma G^k (R(s,s^\prime)) f(s^\prime) d s^\prime\,,\\
    \op D^k f(s) &\coloneqq \text{p.v.} \int_\gamma \frac{\partial}{\partial  n(s^\prime)} G^k (R(s,s^\prime)) f(s^\prime) d s^\prime\,,\\
    \op D^{*k} f(s) &\coloneqq \text{p.v.} \int_\gamma \frac{\partial}{\partial  n (s)} G^k (R(s,s^\prime)) f(s^\prime) d s^\prime\,,
    \end{align}
    \begin{equation}
    \op N^k f(s) \coloneqq - \frac{\partial}{\partial  n (s)} \int_\gamma \frac{\partial}{\partial  n(s^\prime)} G^k (R(s,s^\prime)) f(s^\prime) d s^\prime\,,
\end{equation}
where $R\left( s , s' \right)=\sqrt{\left(\gamma(s')-\gamma(s)\right)^\T\cdot\left(\gamma(s')-\gamma(s)\right)}$.
The two-dimensional free-space Green's function is $G^k(R) = \frac{\ii}{4}H_0^{(1)}(kR)$, where $H_0^{(1)}$ is the Hankel function of the first kind \cite{abramowitz1964handbook}.
These are the building blocks for the electric and magnetic field integral equations (EFIE and MFIE), that relate the longitudinal and transverse electric current densities $J_z$ and $J_t$ to the impinging electromagnetic fields $(E_z,H_t)$ and $(E_t,H_z)$ \cite{jackson1999classical}.

By combining and preconditioning the EFIE and MFIE, we form the Calder\'{o}n combined field integral equation \cite{andriulli2015high}, hereafter denoted ``CCFIE'', which reads for TM and TE polarizations respectively,
\begin{align}
    &\left[{\frac{k}{\tilde{k}}\op N^{\tilde{k}}} \op S^k  + \left(\frac{1}{2} \op I -  {\op D}^{*\tilde{k}} \right)\left(\frac{1}{2} \op I +  \op D^{*k} \right)\right] (J_z) (s) \nonumber\\&\quad= -\frac{\op N^{\tilde{k}}}{\tilde{k}}\frac{E_z(s)}{\ii\eta}+\left(\frac{1}{2} \op I -  {\op D}^{*\tilde{k}} \right)H_t(s)\,,\label{eqn:TMCCFIE}\\
    &\left[{\frac{\tilde{k}}{k}\op S^{\tilde{k}}} \op N^k  + \left(\frac{1}{2} \op I +  {\op D^{\tilde{k}}} \right)\left(\frac{1}{2} \op I -  \op D^k \right)\right] (J_t) (s) \nonumber\\&\quad=\tilde{k}{\op S^{\tilde{k}}}\frac{E_t(s)}{\ii\eta}-\left(\frac{1}{2} \op I +  {\op D^{\tilde{k}}} \right)H_z(s)\,,\label{eqn:TECCFIE}
\end{align}
where $\tilde{k}(s) = k+\ii k_i(s)$ with $k_i(s) \ll k$ for $k \rightarrow \infty$; when $\gamma$ is the circle of radius $a$, a convenient choice is $k_i = 24^{-1/3}k^{1/3}a^{-2/3}$ \cite{antoine2006improved}.

Given a plane wave traveling along direction $p$ and impinging  on $\gamma$, the expressions for the incident fields are given by
\begin{align}
    E_z(s) &= E_0 e^{\ii k p \cdot \gamma(s)}\\
    H_t(s) &= -\frac{E_0}{\eta} e^{\ii k p \cdot \gamma(s)}  \left( p \cdot n(s)\right)
\end{align}
when the polarization is TM, and by
\begin{align}
    E_t(s) &= \ii  E_0  e^{\ii k p \cdot \gamma(s)} \left( p \cdot n(s)\right)\\
    H_z(s) &= \ii \frac{E_0}{\eta} e^{\ii k p \cdot \gamma(s)}
\end{align}
when the polarization is TE.

Finally, given the integral operator $\mathcal{Z}$
\begin{equation}
\op Z f(s) = \int_\gamma \text{ker} (R(s,s^\prime)) f(s^\prime) d s^\prime
\end{equation}
and the spectral variable $\xi$, related to the spectral index $q$ as $\xi = 2\pi q/L$ for $ q \in \mathbb{Z}$,
we define the local operator symbol as
\begin{equation}
\sigma^{\op Z}\left( \xi , s \right)=\int_\gamma
\text{ker}\left(R\left( s , s' \right)\right)
\, e^{ i \xi \left( s' - s \right) }\, ds' \,,
\label{eqn:symboldef}
\end{equation}
so that $\mathcal{Z}\left(e^{ i \xi s }\right)=\sigma^{\op Z}\left( \xi , s \right)\, e^{ i \xi s }$.

\section{Microlocal Analysis Away from Glancing}

The semiclassical principal symbols of the integral operators defined above can be obtained by integrating along $\mathbb{R}$ the zeroth-order approximation of their kernels multiplied by the canonical exponential,
\begin{align}
\sigma_\NGl^{\op S^k}(\xi)&=\int_{\mathbb R}\frac{\ii}{4}H_0^{(1)}(k|x|)e^{\ii \xi x}dx=\frac{\ii}{2\sqrt{k^2-\xi^2}} \,,\label{eqn:Sprinc}\\
\sigma_\NGl^{\op D^k}(\xi)&=\sigma_\NGl^{\op D^{*k}}(\xi)=0\,,\\
\sigma_\NGl^{\op N^k}(\xi)&=-\int_{\mathbb{R}}\frac{\ii k}{4}\frac{H_1^{(1)}(k|x|)}{|x|}e^{\ii\xi x}dx =-\frac{\ii}{2}\sqrt{k^2-\xi^2}\,.\label{eqn:Nprinc}
\end{align}

These results can be directly employed in approximating the preconditioned RHSs of equations \eqref{eqn:TMCCFIE} and \eqref{eqn:TECCFIE} away from the shadow boundary given a plane wave excitation. After defining the local tangential frequency $\xi_T(s)$ as the derivative of the plane wave phase, $\xi_T(s)\coloneq k p\cdot \tau(s)$, we obtain that, by substituting the principal symbols of the operators, the RHSs of equations \eqref{eqn:TMCCFIE}, \eqref{eqn:TECCFIE} are approximated by
\begin{align}
\eqref{eqn:TMCCFIE} &\sim -\frac{E_0}{2\eta}e^{\ii k p\cdot\gamma(s)} (p\cdot n(s))
\left[1-\frac{\sqrt{\tilde{k}^2-\xi_T(s)^2}}{\tilde k (p\cdot n(s))}\right]  \,, \\
\eqref{eqn:TECCFIE} &\sim -\ii\frac{E_0}{2\eta}e^{\ii k p\cdot\gamma(s)}
\left[1-\frac{\tilde k (p\cdot n(s))}{\sqrt{\tilde{k}^2-\xi_T(s)^2}}\right]\,.
\end{align}
Let $\delta:=k_i/k$, then $(\tilde k^{\,2}-\xi_T(s)^2)=k^2\left(p\cdot n(s)\right)^2+\ii 2  k k_i-k_i^2$, so that
\begin{align}
\left[1-\frac{\sqrt{\tilde{k}^2-\xi_T(s)^2}}{\tilde k (p\cdot n(s))}\right]
&=
1-\frac{\sqrt{\left(p\cdot n(s)\right)^2+\ii2 \delta-\delta^2}}{(1+\ii\delta)\,\left(p\cdot n(s)\right)}\,,\label{eqn:portTM}\\
\left[1-\frac{\tilde k (p\cdot n(s))}{\sqrt{\tilde{k}^2-\xi_T(s)^2}}\right]
&=
1-\frac{(1+\ii\delta)\,\left(p\cdot n(s)\right)}{\sqrt{\left(p\cdot n(s)\right)^2+\ii 2 \delta-\delta^2}}.\label{eqn:portTE}
\end{align}
For $2\delta \ll \left(p\cdot n(s)\right)^2$, then
\begin{align}
    &\sqrt{\left(p\cdot n(s)\right)^2+\ii 2 \delta-\delta^2} =\nonumber\\&\quad\quad |\left(p\cdot n(s)\right)|\left(1+\ii \frac{\delta}{\left(p\cdot n(s)\right)^2}+\mathcal{O}\left(\frac{\delta^2}{\left(p\cdot n(s)\right)^4}\right)\right)
\end{align}
and
\begin{equation}
    \eqref{eqn:portTM} =\eqref{eqn:portTE} = 2\Pi_L(s,p) + \mathcal{O}\left(\frac{\delta}{\left(p\cdot n(s)\right)^2}\right)\,,
\end{equation}
where $\Pi_L(s,p)$ is the lit gate function
\begin{equation}
    \Pi_L(s,p) = \begin{cases}
			1, & \text{if $p\cdot n(s) \le 0$}\\
            0, & \text{otherwise}
		 \end{cases}\,.
\end{equation}
By choosing $\delta(s) \propto k^{-2/3}\kappa(s)^{2/3}$, as is typical \cite{antoine2006improved,consoli2022fast}, the condition $|p\cdot n(s)| \gg \sqrt{2\delta(s)}$ results in $|p\cdot n(s)| \gg k^{-1/3}\kappa(s)^{1/3}$, which, by expanding $p\cdot n(s)$ around $s_0$ such that $p\cdot n(s_0) = 0$, results in $|s-s_0| \gg k^{-1/3}\kappa(s)^{-2/3}$. Hence, this analysis is valid outside the Fock region of creeping waves \cite{fock1965electromagnetic}, or spatial glancing region, where, in the high-frequency limit, the RHSs of equations \eqref{eqn:TMCCFIE}, \eqref{eqn:TECCFIE} converge, apart from a constant factor, to the physical optics solutions.

Consistently, by substituting the principal symbols of the operators into the left-hand sides (LHSs) of equations \eqref{eqn:TMCCFIE}, \eqref{eqn:TECCFIE}, we recognize that the principal symbols of the Calder\'{o}n combined operators equal $\op I/2+\mathcal{O}(\delta)$ away from the spectral glancing region. Hence, the principal symbol of the CCFIO converges to the scaled identity as long as $\delta \ll 1$ in the high-frequency limit, confirming the assumption \eqref{eqn:CCFIOassumption} in \Cref{sec:introduction}.

\section{Microlocal Analysis at Glancing}
Whenever the spatial oscillation frequency of the integrand is high enough, equation \eqref{eqn:symboldef} can be approximated as
\begin{equation}
\sigma^{\op Z}\left( \xi , s \right)\simeq\int_\mathbb{R}
\text{ker}\left(R\left( s , s' \right)\right)
\, e^{ i \xi \left( s' - s \right) }\, \chi\left(\frac{s'-s}{\varepsilon}\right)\, ds' \,,
\label{eqn:SymbolFourier}
\end{equation}
where $\chi(t)\in C^\infty(\mathbb R)$ is an even gate function equal to $1$ for $|t|<1/2$ and zero for $|t|>1$. By setting $\varepsilon\propto k^{-1/3}\kappa(s)^{-2/3}$ to cover the spatial glancing zone, $kR\rightarrow \infty$ over a significant portion of $\left(s-\varepsilon,s+\varepsilon\right)$ for $k\rightarrow \infty$. By substituting the Hankel function in the integrand kernel with its large argument asymptotic expansion \cite{abramowitz1964handbook}, we obtain that, in the glancing spectral region, i.e., at $\xi \sim k$, the semiclassical transition symbols of the integral operators are approximated as \cite{melrose1975local}
\begin{align}
    \sigma_\Gl^{{\op S}^k}(\xi,s) &\simeq \frac{e^{\ii\pi/4}}{2\sqrt{\pi}}\left(\frac{k^2\kappa(s)}{\sqrt{3}}\right)^{-1/3}F\left(K\right)\label{eqn:Sgl_F}\\
    &= {\ii \pi}\left(2k^2\kappa(s)\right)^{-1/3} \Ai(x)\left(\Ai(x)-\ii \Bi(x)\right)\,,\label{eqn:Sgl_Ai}\\
    \sigma_\Gl^{{\op D}^k}(\xi,s)&=\sigma_\Gl^{{\op D}^{*k}}(\xi,s) \simeq
    \frac{e^{\ii\pi/4}}{2\sqrt{\pi}} \sqrt{3}\, F^\prime(K)
    \label{eqn:Dgl_F}\\
    &= -\ii \frac{\pi}{2} \left[\Ai(x)\left(\Ai(x)-\ii \Bi(x)\right)\right]^\prime\,,
    \label{eqn:Dgl_Ai}\\
    \sigma_\Gl^{{\op N}^k}(\xi,s) &\simeq -\frac{e^{\ii\pi/4}}{2\sqrt{\pi}}\left(\frac{{k^2\kappa(s)}}{\sqrt{3}}\right)^{1/3}\left(KF(K)+6F''(K)\right)\label{eqn:Ngl_F}\\
    &= -\ii \pi \left(2k^2\kappa(s)\right)^{1/3} \Ai'(x)\left(\Ai'(x)-\ii\Bi'(x)\right)\,,\label{eqn:Ngl_Ai}
\end{align}
where $K=K(\xi,s):=(k-\xi)24^{1/3}k^{-1/3}\kappa(s)^{-2/3}$ and $F(K)\coloneq
\int_0^\infty
u^{-1/2}e^{\ii\left(Ku-u^3\right)}\,du$. Moreover, $x=x(\xi,s)\coloneq-{K}/{12^{1/3}}$ and $\Ai(x)$, $\Bi(x)$ are Airy functions \cite{abramowitz1964handbook}.
By substituting the Airy Wronskian \cite{abramowitz1964handbook}, from \eqref{eqn:Dgl_Ai} we further obtain
\begin{align}
    \frac{1}{2}+\sigma_\Gl^{{\op D}^k} &\simeq -\ii \pi \Ai^\prime(x)\left(\Ai(x)-\ii \Bi(x)\right)\,,\\
    \frac{1}{2}-\sigma_\Gl^{{\op D}^k} &\simeq\ii \pi \Ai(x)\left(\Ai(x)-\ii \Bi(x)\right)^\prime\,.
\end{align}
Symmetric results can be obtained for the negative spectral indices branch of eigenvalues, that is, for $\xi \sim -k$.

By applying stationary phase and endpoint approximations to the expressions \eqref{eqn:Sgl_F}, \eqref{eqn:Dgl_F}, \eqref{eqn:Ngl_F}, or, equivalently, by using the large argument asymptotic expansions of the Airy functions \cite{abramowitz1964handbook} in \eqref{eqn:Sgl_Ai}, \eqref{eqn:Dgl_Ai}, \eqref{eqn:Ngl_Ai}, one finds that these symbols can be further approximated as
\begin{equation}
    \sigma_\Gl^{{\op S}^k}(\xi,s) \sim \begin{cases}
        \frac{1}{\sqrt{8k}}(k-\xi)^{-1/2}\left(\ii+  e^{\ii \frac{4}{3}\left(-x\right)^{3/2}}\right), & \text{if $\xi < k$}\\
        \frac{1}{\sqrt{8k}}(\xi-k)^{-1/2}, & \text{if $\xi > k$}
    \end{cases}
    \label{eqn:SglancApprox}
    \end{equation}
\begin{align}
    \sigma_\Gl^{{\op D}^k}(\xi,s) &\sim \begin{cases}
         \frac{\ii}{2}\,e^{\ii \frac{4}{3}\left(-x\right)^{3/2}}, & \text{if $\xi < k$}\\
        \frac{\sqrt{2}}{16}\kappa(s)\sqrt{k}(\xi-k)^{-3/2}, & \text{if $\xi > k$}
    \end{cases}
    \label{eqn:DglancApprox}
    \\
    \sigma_\Gl^{{\op N}^k}(\xi,s)&\sim \begin{cases}
        \frac{\sqrt{k}}{\sqrt{2}}(k-\xi)^{1/2}\left(-\ii+e^{\ii \frac{4}{3}\left(-x\right)^{3/2}}\right), & \text{if $\xi < k$}\\
        \frac{\sqrt{k}}{\sqrt{2}}(\xi-k)^{1/2}, & \text{if $\xi > k$}
    \end{cases}\,.
    \label{eqn:NglancApprox}
\end{align}

In the high-frequency limit, these approximations are uniformly valid for constant values of $K$, that is, for $|k-\xi|\sim 24^{-1/3}k^{1/3}\kappa(s)^{2/3}$, consistently with an increase of the spectral glancing region in frequency as $k^{1/3}\kappa(s)^{2/3}$.

After denoting
\begin{align}
    \tilde{x}&\coloneq (\xi-\tilde{k})2^{1/3}\tilde{k}^{-1/3}\kappa(s)^{-2/3}\nonumber\\
    &=\left(\frac{k}{\tilde{k}}\right)^{1/3}\left(x-\ii\frac{k_i}{2^{-1/3}k^{1/3}\kappa(s)^{2/3}}\right)\,,
\end{align}
we note that the sum of products of the glancing symbols of the operators forming the CCFIO is asymptotically constant at $x = 0$ in frequency and curvature if $k_i \propto k^{1/3}\kappa(s)^{2/3}$ for $k \rightarrow \infty$, corresponding to constant values of $\tilde{x}$. Moreover, for $x$ bounded in an asymptotically constant interval symmetrically centered around $0$, corresponding to the spectral glancing region, $\tilde{x}$ is also bounded in a constant subset of $\mathbb{C}$ if $k_i \propto k^{1/3}\kappa(s)^{2/3}$ for $k \rightarrow \infty$. Hence, we infer that, in this case, the dominant part of the spectral glancing symbol of the CCFIO is asymptotically bounded in the high-frequency regime, due to the boundedness of the variables $x$ and $\tilde{x}$, resulting in a frequency increase of the spectral content of the operator $\left(\text{CCFIO}-\op I /2\right)$ at least proportional to the size of the glancing zone. As the CCFIO converges, in principal value, to the scaled identity outside the spectral glancing region, we finally conclude that the spectral information of the compact part  $\left(\text{CCFIO}-\op I /2\right)$ is limited to the spectral glancing region, increasing as $k^{1/3}$ in the high-frequency limit, which corroborates the estimate in \Cref{sec:introduction}.

Notice that \eqref{eqn:SymbolFourier} represents a windowed, or short-time, Fourier transform of the kernel $\text{ker}(R)$. By multiplying the inverse of the glancing symbol $\sigma_\Gl^{{\op S}^k}(\xi,s)$ by the windowed Fourier transform of the TM-EFIE RHS and inverse transforming the result, we obtain a valid approximation of the TM current $J_z$ at glancing,
\begin{align}
    J_{z,\Gl}(s_0+t)&\sim \frac{2E_0}{\eta}\,\left(2 k^{-1}\kappa(s_0)\right)^{1/3}\,
e^{ik(p \cdot \gamma(s_0))}\,
e^{ik t}\,\nonumber\\
&\quad
I^{\text{TM}}\!\left(\left(\frac{k \kappa(s_0)^2}{2}\right)^{1/3}t\right)
\label{eqn:Jzapprox}
\end{align}
and, similarly from $\sigma_\Gl^{{\op N}^k}(\xi,s)$, a valid approximation of the TE current $J_t$ at glancing,
\begin{align}
    J_{t,\Gl}(s_0+t)&\sim \frac{2E_0}{\eta}\,
e^{ik(p \cdot \gamma(s_0))}\,
e^{ik t}\,\nonumber\\
&\quad
I^{\text{TE}}\!\left(\left(\frac{k \kappa(s_0)^2}{2}\right)^{1/3}t\right)\,,
\label{eqn:Jtapprox}
\end{align}
where $s_0$ is such that $p\cdot n(s_0)=0$ and $p\cdot \tau(s_0)>0$, that is, $s_0$ is in the Fock glancing region associated to the branch of eigenvalues with positive indices, $|t|<\mathcal{O}(k^{-1/3}\kappa(s_0)^{-2/3})$, and
\begin{align}
    I^{\text{TM}}(\psi)&\coloneq\frac{1}{2\pi}\int_{\mathbb R}
\frac{e^{i\psi x}}{\Ai(x)-\ii\Bi(x)}\,dx\label{eqn:ITM_approx}\\
I^{\text{TE}}(\psi)&\coloneq\frac{1}{2\pi}\int_{\mathbb R}
\frac{e^{i\psi x}}{\Ai^\prime(x)-\ii\Bi^\prime(x)}\,dx\,.\label{eqn:ITE_approx}
\end{align}

Note that this approximate procedure of \textit{solving by symbol} can be seen as a generalization of the exact procedure of \textit{solving by eigenvalues} valid only for circular geometries.

\section{Numerical Results}

\begin{figure}
\subfloat[\label{subfig-1:}]{%
  \includegraphics[width=0.9\columnwidth]{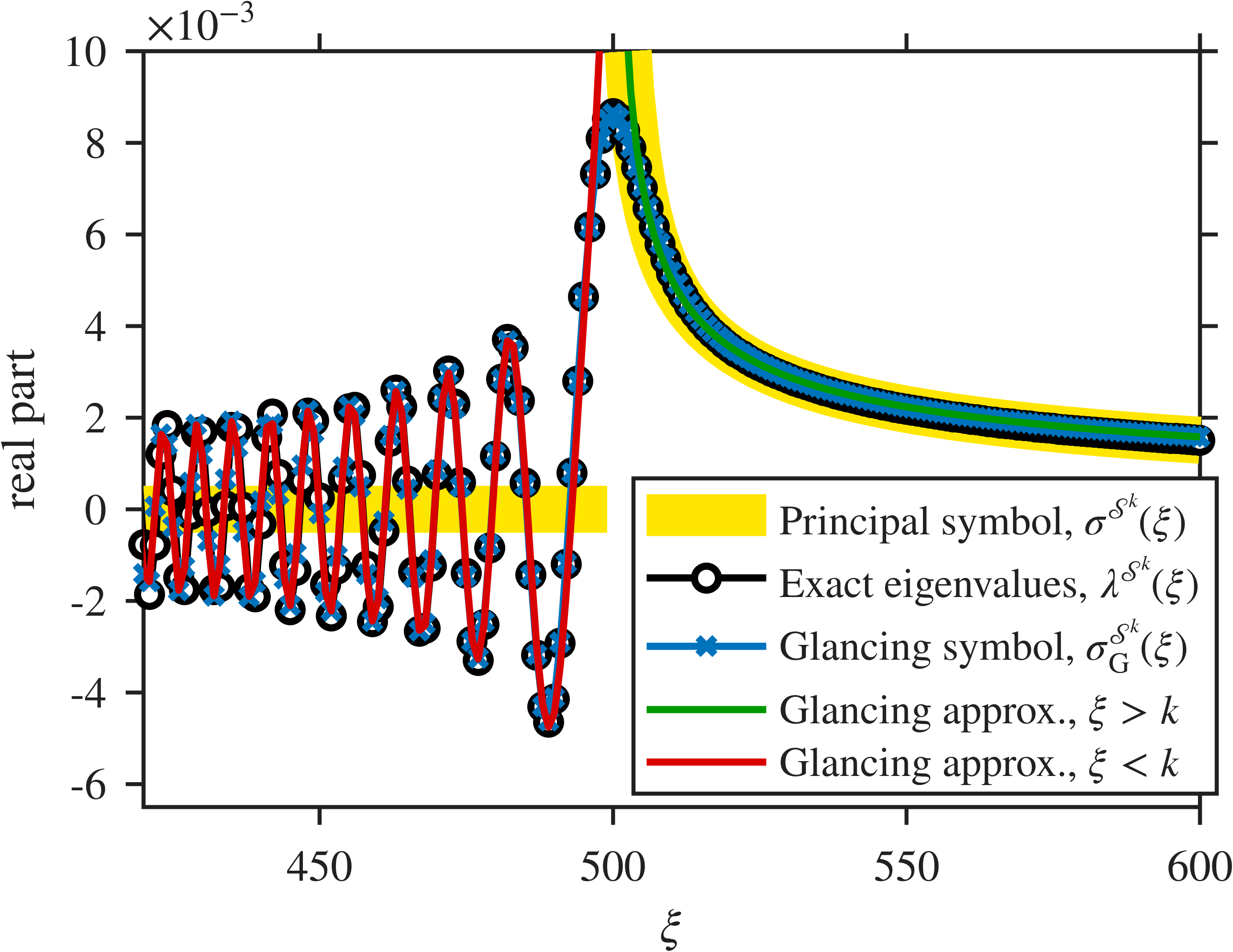}
}
\hfill
\subfloat[\label{subfig-2:}]{%
  \includegraphics[width=0.9\columnwidth]{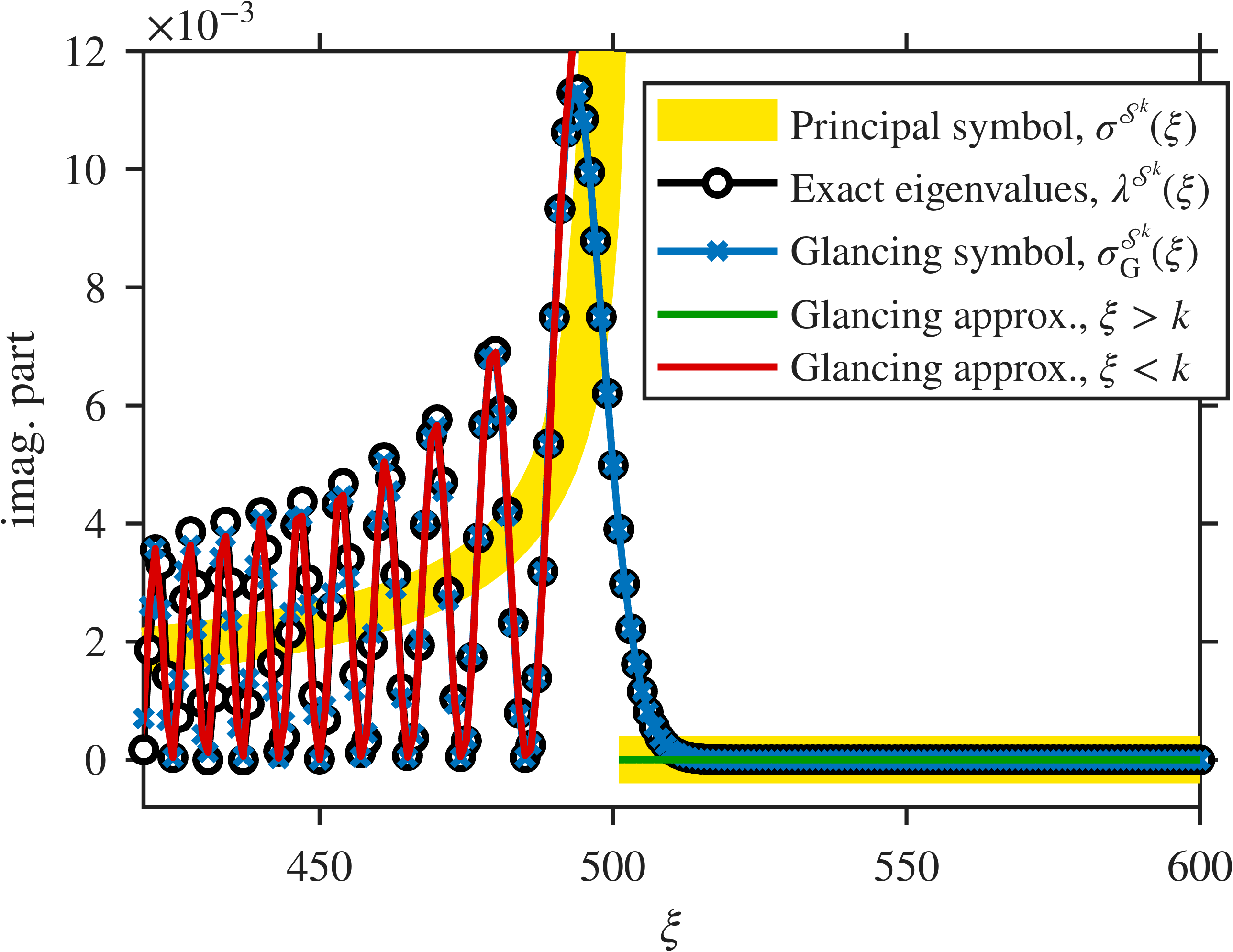}
}
\caption{Comparison between the principal symbol \eqref{eqn:Sprinc}, the glancing symbol \eqref{eqn:Sgl_Ai}, the approximate glancing symbol \eqref{eqn:SglancApprox} and the exact eigenvalues $\lambda^{\op S^k}(\xi) = \ii \pi a /2 \,J_{a\xi}(ka)H_{a\xi}^{(1)}(ka)$ of the single-layer operator evaluated over the circular boundary at $ka = 500$: real part (top) and imaginary part (bottom).}
\label{fig:spectrumSL}
\end{figure}

\begin{figure}
\subfloat[\label{subfig-1:}]{%
  \includegraphics[width=0.9\columnwidth]{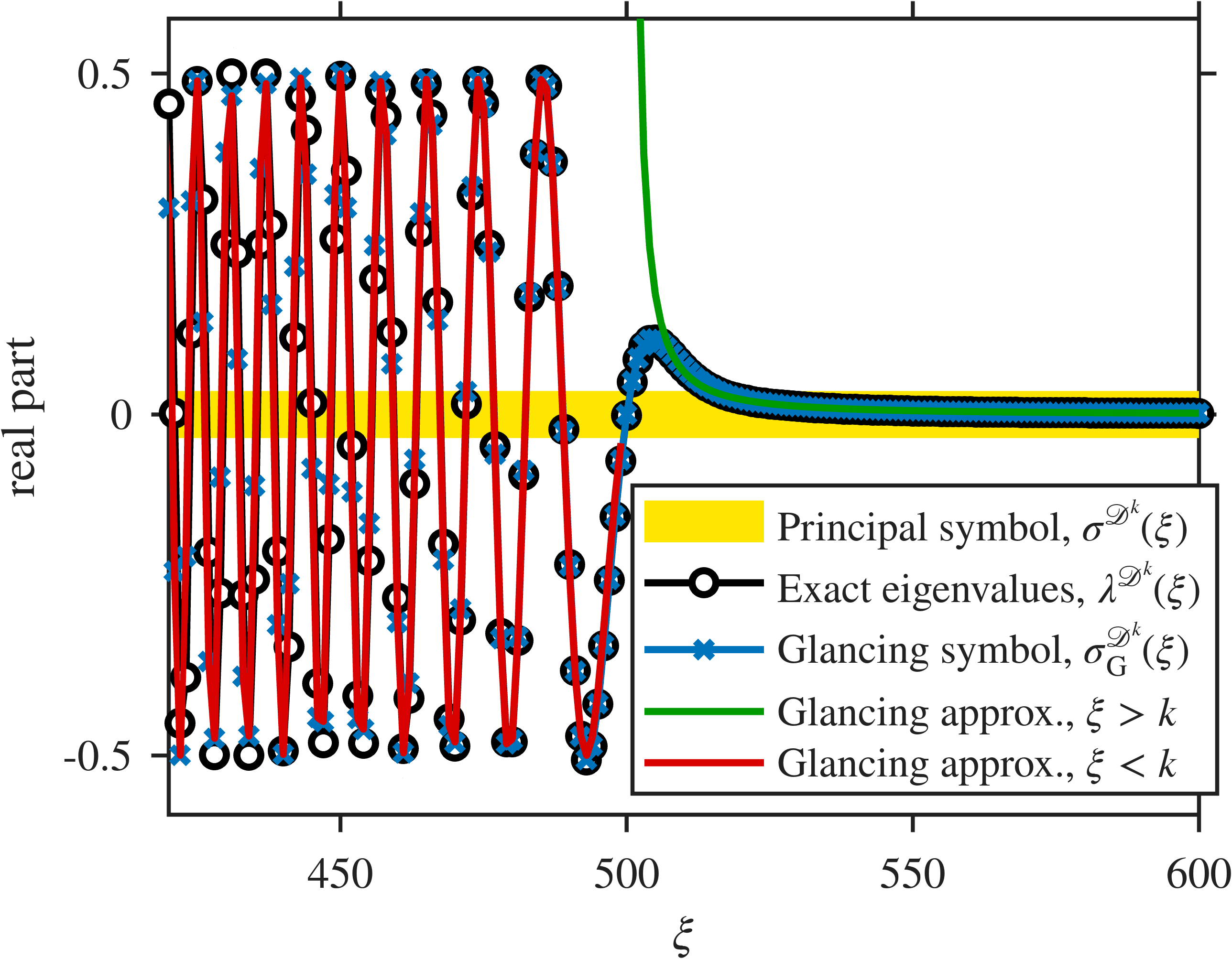}
}
\hfill
\subfloat[\label{subfig-2:}]{%
  \includegraphics[width=0.9\columnwidth]{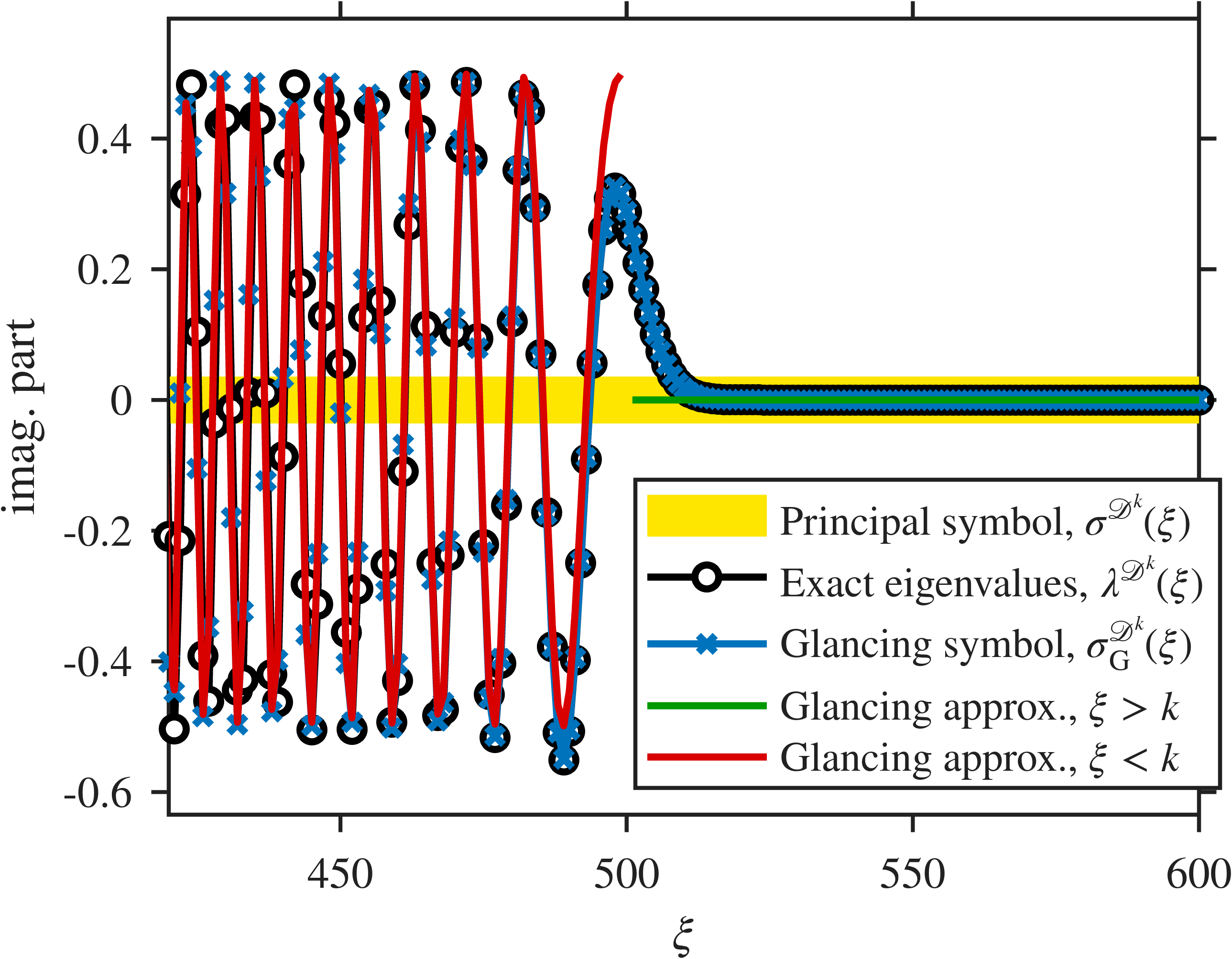}
}
\caption{Comparison between the glancing symbol \eqref{eqn:Dgl_Ai}, the approximate glancing symbol \eqref{eqn:DglancApprox} and the exact eigenvalues $\lambda^{\op D^k}(\xi) = \ii \pi a k/4 \,\left(J_{a\xi}(ka)H_{a\xi}^{(1)}(ka)\right)^\prime$ of the double-layer operator evaluated over the circular boundary at $ka = 500$: real part (top) and imaginary part (bottom).}
    \label{fig:spectrumDL}
\end{figure}

\begin{figure}
\subfloat[\label{subfig-1:}]{%
  \includegraphics[width=0.9\columnwidth]{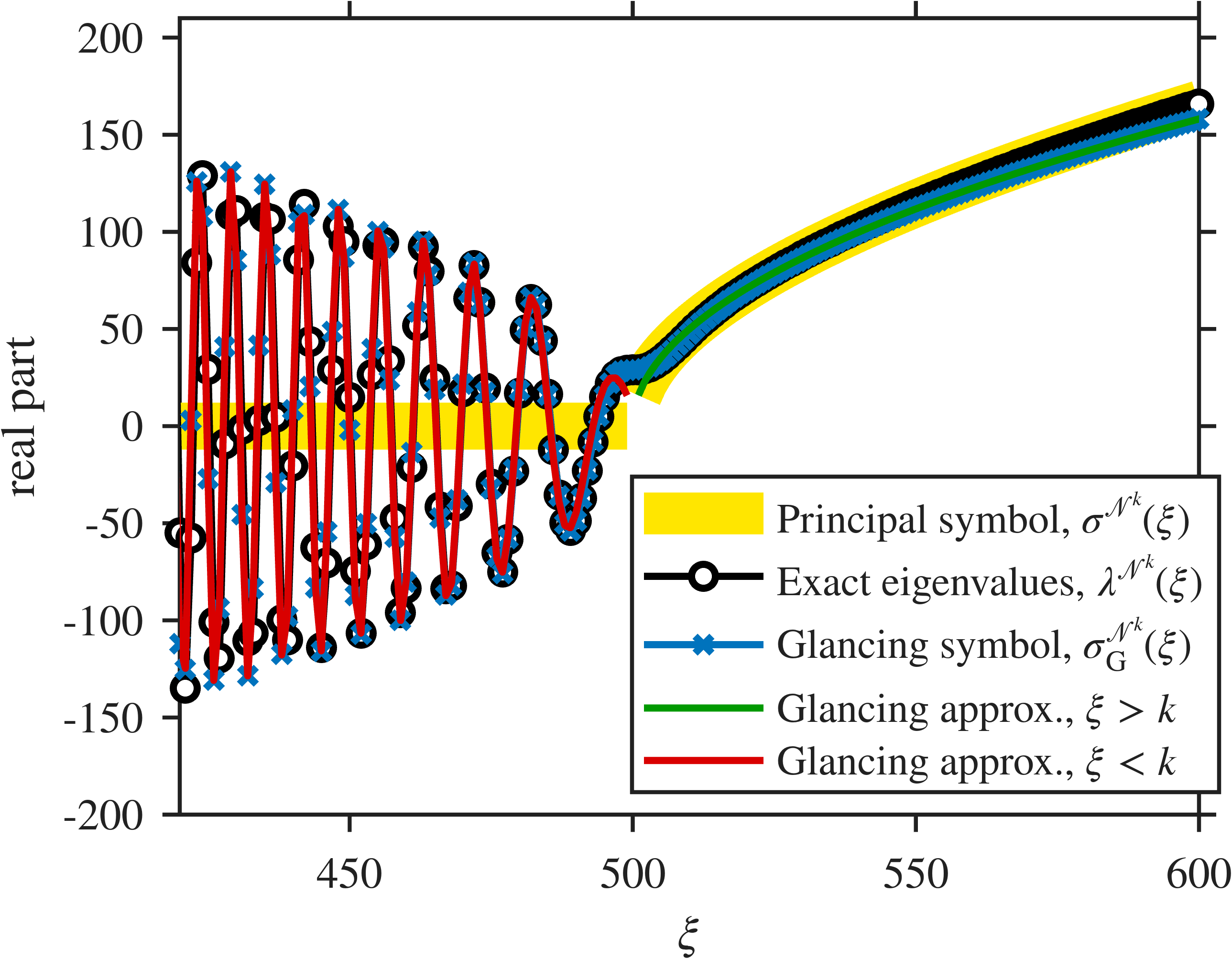}
}
\hfill
\subfloat[\label{subfig-2:}]{%
  \includegraphics[width=0.9\columnwidth]{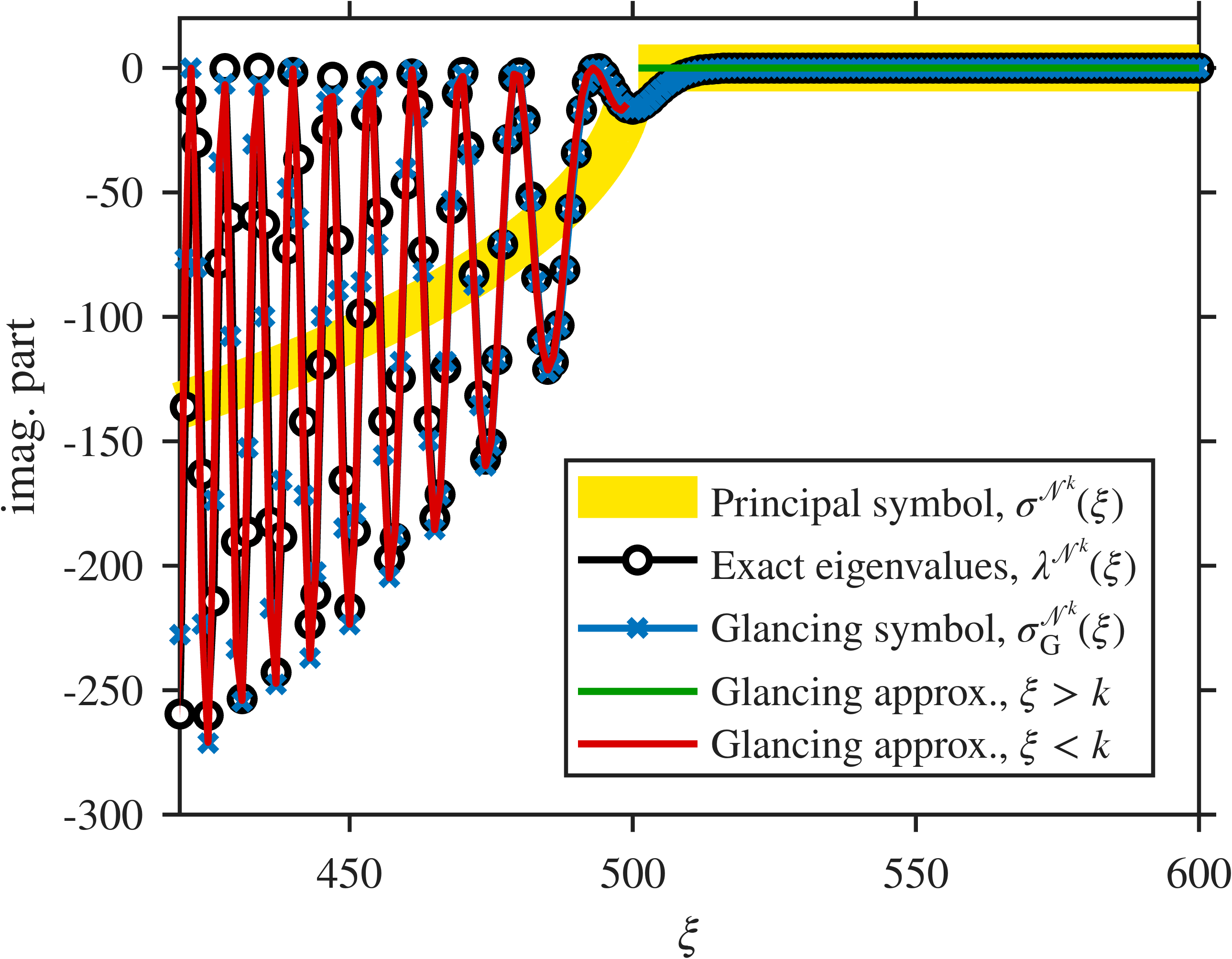}
}
\caption{Comparison between the principal symbol \eqref{eqn:Nprinc}, the glancing symbol \eqref{eqn:Ngl_Ai}, the approximate glancing symbol \eqref{eqn:NglancApprox} and the exact eigenvalues $\lambda^{\op N^k}(\xi) = -\ii \pi a  k^2/2 \,J^\prime_{a\xi}(ka)H_{a\xi}^{(1)\prime}(ka)$ of the hypersingular operator evaluated over the circular boundary at $ka = 500$: real part (top) and imaginary part (bottom).}
    \label{fig:spectrumHS}
\end{figure}

The numerical results aim to validate the general framework introduced in the contribution. First, we compare the principal and glancing symbols with the exact eigenvalue spectra of the single-layer operator (\Cref{fig:spectrumSL}), double-layer operator (\Cref{fig:spectrumDL}), and hypersingular operator (\Cref{fig:spectrumHS}) applied over the circular boundary. In all the cases, we observe good correspondence between the exact eigenvalues and the glancing symbols in the glancing spectral region (i.e., around $\xi = k$); moreover, a good agreement is found even with the approximations in \eqref{eqn:SglancApprox}, \eqref{eqn:DglancApprox}, \eqref{eqn:NglancApprox}.

\begin{figure}
\subfloat[\label{subfig-1:}]{%
  \includegraphics[width=0.9\columnwidth]{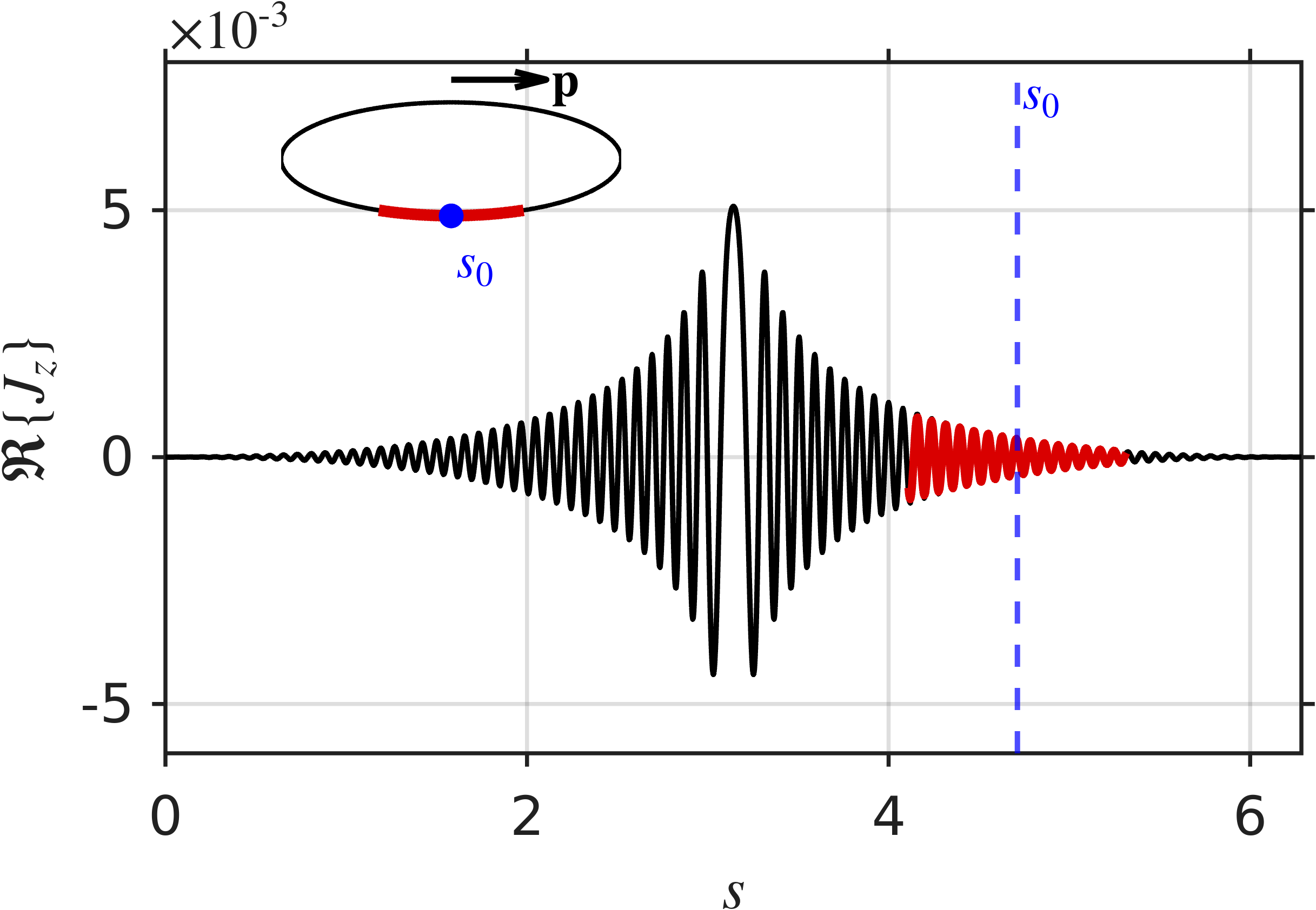}
}
\hfill
\subfloat[\label{subfig-2:}]{%
  \includegraphics[width=0.9\columnwidth]{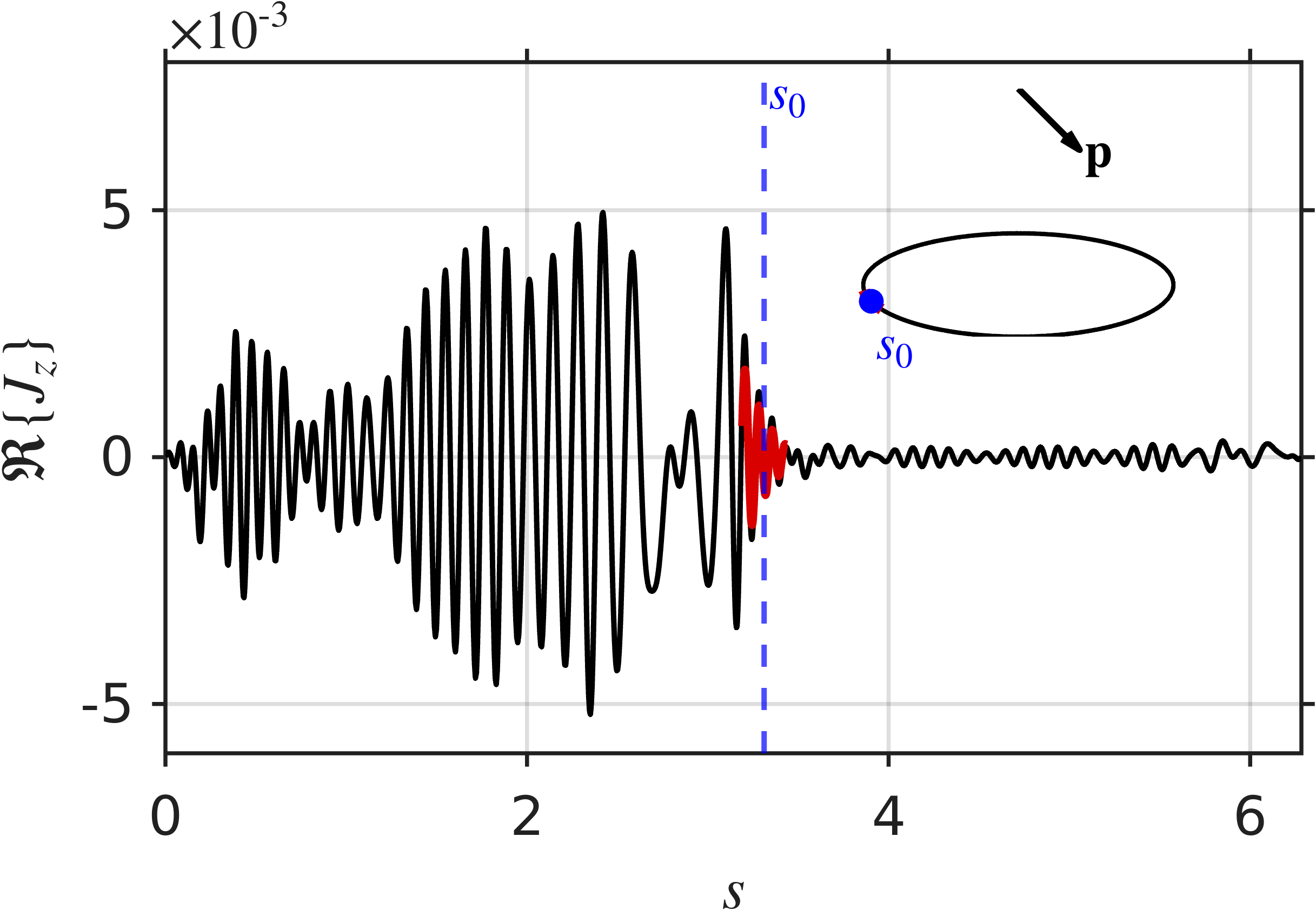}
}
\hfill
\subfloat[\label{subfig-2:}]{%
  \includegraphics[width=0.9\columnwidth]{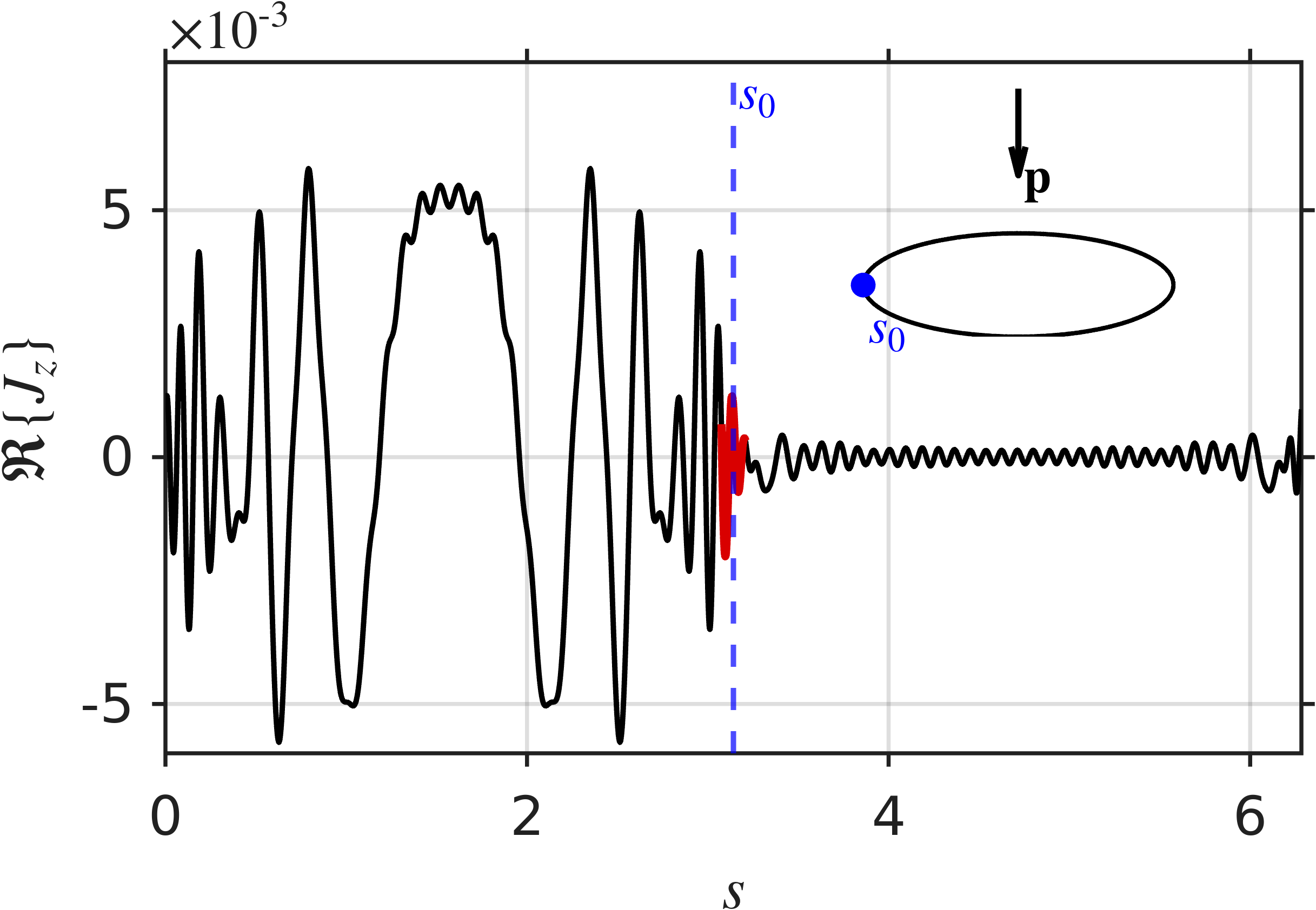}
}
\caption{Comparison between the exact current (in black) and its approximation around $s_0$ \eqref{eqn:Jzapprox} (in red) for different plane wave incidence directions ($kL/(2\pi)=80$).}
    \label{fig:TMglancingcurrent}
\end{figure}

\begin{figure}
\subfloat[\label{subfig-1:}]{%
  \includegraphics[width=0.9\columnwidth]{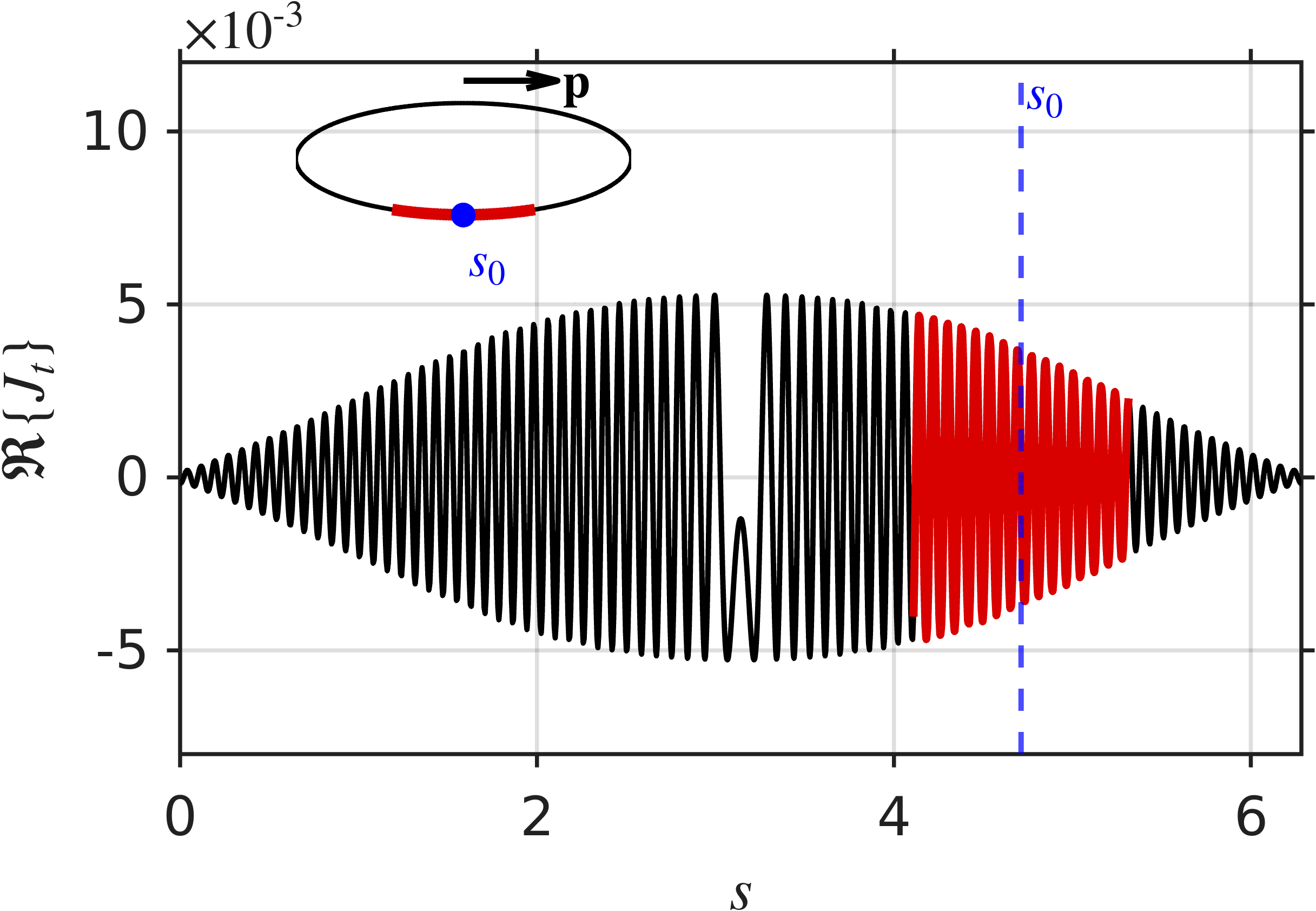}
}
\hfill
\subfloat[\label{subfig-2:}]{%
  \includegraphics[width=0.9\columnwidth]{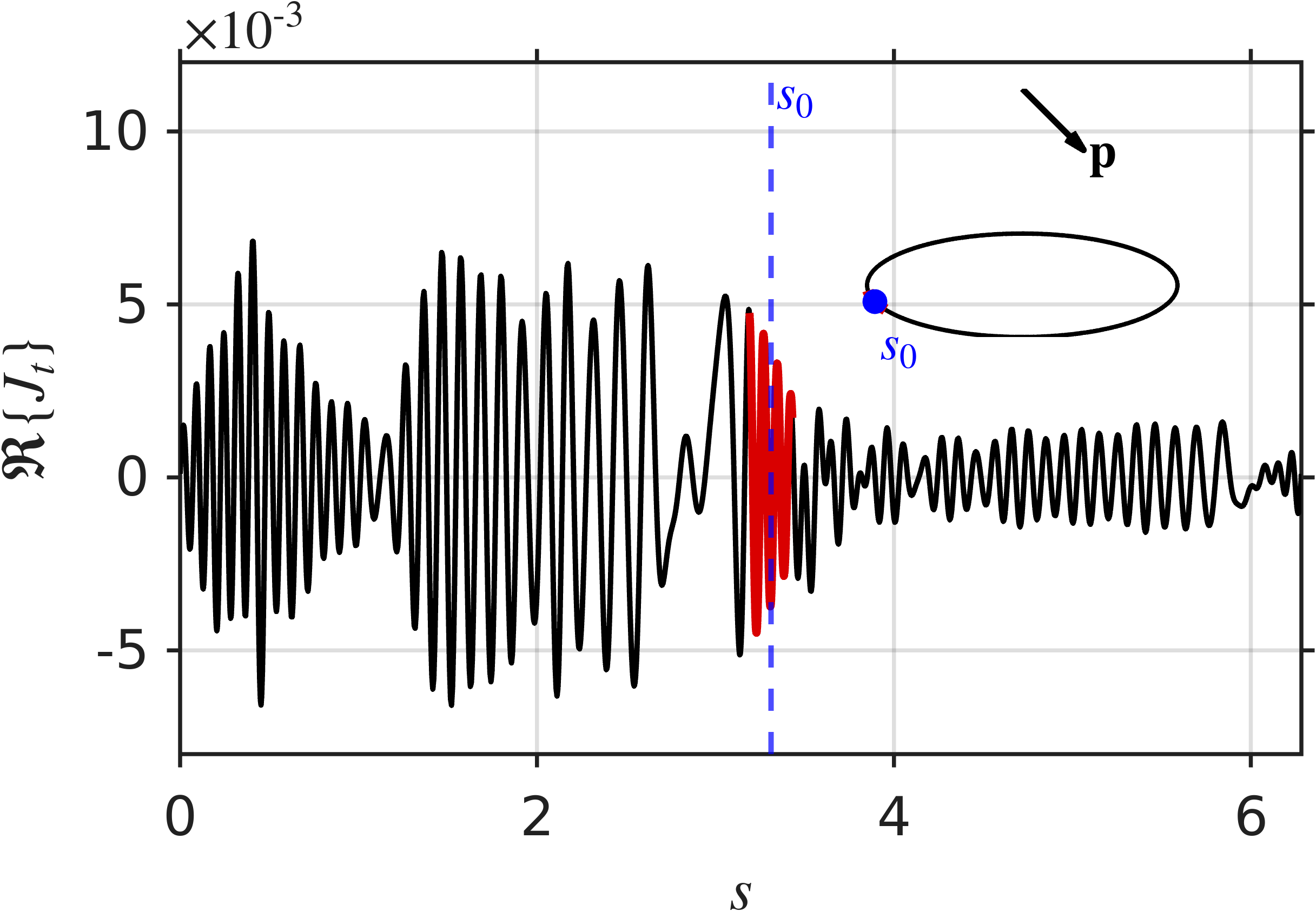}
}
\hfill
\subfloat[\label{subfig-2:}]{%
  \includegraphics[width=0.9\columnwidth]{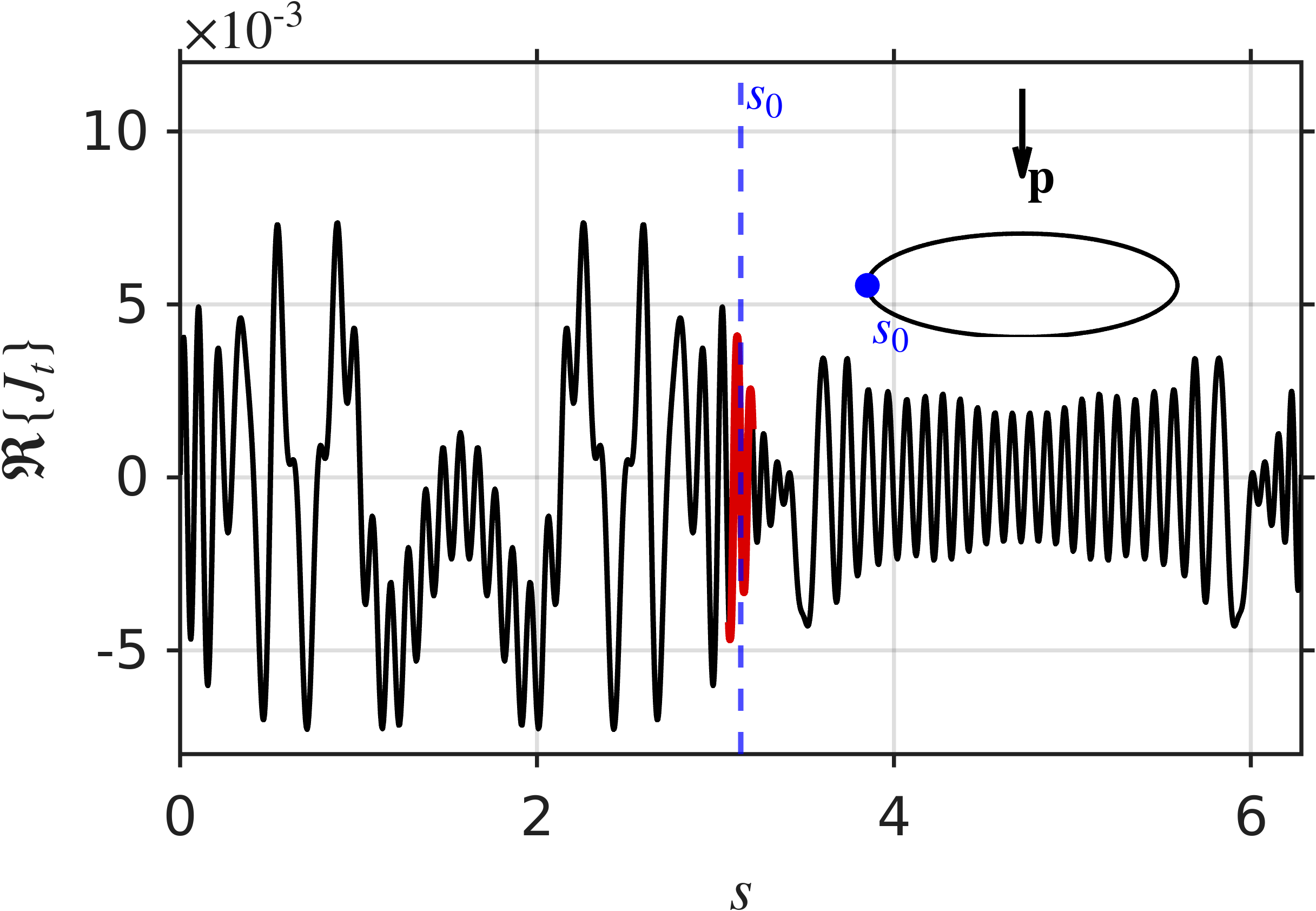}
}
\caption{Comparison between the exact current (in black) and its approximation around $s_0$ \eqref{eqn:Jtapprox} (in red) for different plane wave incidence directions ($kL/(2\pi)=80$).}
    \label{fig:TEglancingcurrent}
\end{figure}

Second, we test the approximation of the glancing currents \eqref{eqn:ITM_approx} and \eqref{eqn:ITE_approx} given a plane wave impinging on an elliptical cylinder (\Cref{fig:TMglancingcurrent} and \Cref{fig:TEglancingcurrent}). We consider different incident directions $p$, resulting in glancing regions of different sizes, given the non-constant curvature of the boundary $\gamma$. In all the cases, the approximate current \eqref{eqn:ITM_approx} (in red) is in good agreement with the exact current (in black), obtained from analytic formulae, in the Fock, or glancing, spatial region of curvilinear abscissas $s$ such that $|s-s_0|\le k^{-1/3}\kappa(s_0)^{-2/3}$, where $s_0$ is such that $p\cdot n(s_0)=0$ and $p\cdot \tau(s_0)>0$.

Overall, these results indicate that the semiclassical microlocal framework described in this contribution provides a useful tool to interpret and predict how the spectral filter acts and impacts the solution accuracy of the fast direct solver as the frequency increases.

\section*{Acknowledgements}
This work was supported by the European Innovation Council (EIC) through the European Union’s Horizon Europe research Programme under Grant 101046748 (Project CEREBRO).

{ \printbibliography}

\end{document}